\documentclass[prl,twocolumn,floatfix]{revtex4}

\usepackage{graphicx}
\usepackage{dcolumn}
\usepackage{amsmath}
\usepackage{amssymb}
\usepackage{bm}
\usepackage[english]{babel}
\usepackage{color}

\renewcommand{\P}{\mathbb{P}}
\newcommand{\Nb}{N_\beta}
\newcommand{\bq}{\begin{equation}}
\newcommand{\eq}{\end{equation}}
\begin{document}

\title{Full characterization of the fractional Poisson process}

\author{Mauro Politi}
\author{Taisei Kaizoji}
\affiliation{ SSRI \& Department of Economics and Business, \\ International Christian University, 3-10-2 Osawa, Mitaka, Tokyo, 181-8585 Japan.}

\author{Enrico Scalas}
\affiliation{Dipartimento di Scienze e Tecnologie Avanzate,\\Universit\`a del Piemonte Orientale ``Amedeo Avogadro'', 
Viale T. Michel 11, 15121 Alessandria, Italy.}

\date{\today}%

\begin{abstract}
The fractional Poisson process (FPP) is a counting process with independent and identically distributed inter-event times following
the Mittag-Leffler distribution. This process is very useful in several fields of applied and theoretical physics including models for
anomalous diffusion. Contrary to the well-known Poisson 
process, the fractional Poisson process does not have stationary and independent increments. It is not a L\'evy process and it is
not a Markov process. In this letter, we present formulae for its finite-dimensional distribution functions, fully characterizing
the process. These exact analytical results are compared to Monte Carlo simulations. 
\end{abstract}


\maketitle

%%%%%%%%%%%%%%%%%%%%%%%%%%%%%%%%%%%%%%%%%%%%%%%%%%%%%%%%%%%%%%%%%%%%

% Introduction
From a loose mathematical point of view, counting processes $N(A)$ are stochastic processes that count the random number of points in a set $A$.
They are used in many fields of physics and other applied sciences. In this letter, we will consider one-dimensional real sets with the 
physical meaning of time intervals. The points will be incoming events whose duration is much smaller than the inter-event or
inter-arrival waiting time. For instance, counts from a Geiger-M\"uller counter can be described in this way. The number of counts, $N(\Delta t)$, 
in a given time interval $\Delta t$ is known to follow the Poisson distribution
\bq
\label{poissonprocess}
\P(N(\Delta t) = n) = \exp(-\lambda \Delta t) \frac{(\lambda \Delta t)^n}{n!},
\eq
where $\lambda$ is the constant rate of arrival of ionizing particles. Together with the assumption of independent and stationary increments, 
Eq. \eqref{poissonprocess} is sufficient to define the {\em homogeneous} Poisson process. Curiously, one of the first occurrences of this process in the
scientific literature was connected to the number of casualties by horse kicks in the Prussian army cavalry \cite{Bortkiewicz1898}. The Poisson
process is strictly related to the exponential distribution. The inter-arrival times $\tau_i$ identically follow the exponential distribution and are independent random variables. This means that the Poisson process is a prototypical {\em renewal} process.
A justification for the ubiquity of the Poisson process has to do with its relationship with the binomial distribution. Suppose that the
time interval of interest $(t,t+\Delta t)$ is divided into $n$ equally spaced sub-intervals. Further assume that a counting event appears in such
a sub-interval with probability $p$ and does not appear with probability $1-p$. Then, $\P(N(\Delta t) = k) = \mathrm{Bin}(k;p,n)$ is a binomial distribution of parameters $p$ and $n$ and the expected number of events in the time interval is
given by $\mathbb{E}[N(\Delta t)] = n p$. If this expected number is kept constant for $n \to \infty$, the binomial distribution converges to the
Poisson distribution of parameter $\lambda = \mathbb{E}[N(\Delta t)]/\Delta t$, while, in the meantime, $p \to 0$. However, it can be shown that
many counting processes with non-stationary increments converge to the Poisson process after a transient period. It is sufficient to require that
they are renewal process (i.e. they have independent and identically distributed (iid) inter-arrival times) and that $\mathbb{E}(\tau_i) < \infty$.
In other words, many counting processes with non-independent and non-stationary increments behave as the Poisson process if observed
long after the transient period. 

% Counting processes with heavy tails
In recent times, it has been shown that heavy-tailed distributed inter-arrival times (for which $\mathbb{E}(\tau_i) = \infty$) do play a role in many phenomena such as blinking nano-dots \cite{Margolin2005,Margolin2006}, 
human dynamics \cite{Barabasi2005,Barabasi2010} and the related inter-trade times in financial markets \cite{Scalas2004b,Scalas2006}. 

\begin{figure}
\includegraphics[width=\columnwidth]{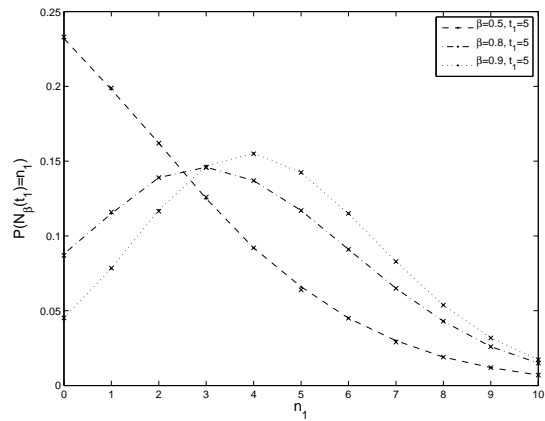}
\caption{$P(\Nb(T_1)=n_1)$ as function of $n_1$ for three different values of $\beta$. The crosses are estimations obtained 
from $10^5$ Monte Carlo samples and the lines are given to guide the eye.}
\label{fig:prob1}
\end{figure}

% Definition of the fractional compound process
Among the counting processes with non-stationary increments, the so-called {\em fractional Poisson process} \cite{Laskin2003}, $\Nb (t)$,
is particularly important because it is the thinning limit of counting processes related to renewal processes with 
power-law distributed inter-arrival times
\cite{Scalas2004,Mainardi2004}. Moreover, it can be used to approximate anomalous diffusion ruled by space-time fractional
diffusion equations \cite{Scalas2004,Metzler2004,Heinsalu2006,Magdziarz2007,Fulger2008,Germano2009,Nane2010}.
It is a straightforward generalization of the Poisson process defined as follows.
Let $\{\tau_i\}_{i=1}^{\infty}$ be a sequence of independent and identically distributed positive random variables with
the meaning of inter-arrival times and let their common cumulative distribution function (cdf) be
\begin{equation}
\label{mlcdf}
F_\tau (t) = \mathbb{P}(\tau \leq t) = 1 - E_\beta (-t^\beta),
\end{equation}
where $E_\beta (-t^\beta)$ is the one-parameter Mittag-Leffler function, $E_\beta(z)$, defined in the complex plane as
\begin{equation}
\label{mlf}
E_\beta (z) = \sum_{n=0}^{\infty} \frac{z^n}{\Gamma(n \beta + 1)}
\end{equation} 
evaluated in the point $z = -t^\beta$ and with the prescription $0 < \beta \leq 1$. In equation \eqref{mlf}, $\Gamma (\cdot)$
is Euler's Gamma function. The sequence of the {\em epochs}, $\{T_n\}_{n=1}^\infty$, is given by the sums of the inter-arrival
times
\begin{equation}
\label{epoch}
T_n = \sum_{i=1}^n \tau_i.
\end{equation}
The epochs represent the times in which events arrive or occur. Let $f_\tau (t) = dF_\tau (t) / dt$ denote the probability
density function (pdf) of the inter-arrival times, then the probability density function of the $n$-th epoch is simply given
by the $n$-fold convolution of $f_\tau (t)$, written as $f^{*n}_\tau (t)$. In Ref. \cite{Mainardi2004}, it is shown that
\begin{equation}
\label{epochpdf}
f_{T_n} (t) = f^{*n}_\tau (t) = \beta \frac{t^{n\beta -1}}{(n-1)!} E^{(n)}_\beta (-t^\beta),
\end{equation}
where $E^{(n)}_\beta (-t^\beta)$ is the $n$-th derivative of $E_\beta (z)$ evaluated in $z=-t^\beta$.
The counting process $\Nb (t)$ counts the number of epochs (events) up to time $t$, assuming that
$T_0 = 0$ is an epoch as well, or, in other words, that the process begins from a renewal point. This
assumption will be used all over this paper. $\Nb (t)$ is given by
\begin{equation}
\label{countingprocess}
\Nb (t) = \max \{ n:\; T_n \leq t \}.
\end{equation}
In Ref. \cite{Scalas2004}, the {fractional Poisson distribution} is derived and it is given by
\begin{equation}
\label{eq:fractionalpoissondistribution}
\P(\Nb (t) = n) = \frac{t^{\beta n}}{n!} E_\beta^{(n)} (-t^\beta).
\end{equation}
Eq.~\eqref{eq:fractionalpoissondistribution} coincides with the Poisson distribution of
parameter $\lambda = 1$ for $\beta =1$.
In principle, equations \eqref{mlf} and \eqref{eq:fractionalpoissondistribution} can be directly used to derive the fractional Poisson
distribution, but convergence of the series is slow. Fortunately, in a recent paper, Beghin and Orsingher proved that
\begin{multline}
\label{beghinorsingher}
E_\beta^{(n)} (-t^\beta) = \\
\frac{n!}{t^{\beta n}} \int_{0}^{\infty} F_{S_\beta} (t;u) \left[ \frac{ \exp(-u) u^{n-1} }{(n-1)!} 
- \frac{\exp(-u) u^n}{n!} \right] \, du,
\end{multline}
where $F_{S_\beta} (t;u)$ is the cdf of a stable random variable $S_\beta(\nu,\gamma,\delta)$  with index $\beta$, 
skewness parameter $\nu = 1$, scale parameter $\gamma=(u \cos\pi\beta/2)^{1/\beta}$ and location $\delta=0$ \cite{Begin2009}. The integral
in equation \eqref{beghinorsingher} can be evaluated numerically and 
Fig. \ref{fig:prob1}1 shows $\mathbb{P}(\Nb(t)=n)$ for three different values of $\beta$. The Monte Carlo
simulation of the fractional Poisson process is based on the algorithm presented in equation (20) of Ref. \cite{Fulger2008}.
\begin{figure}
\includegraphics[width=\columnwidth]{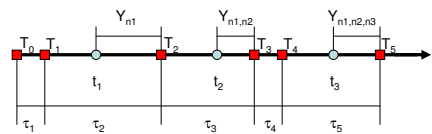}
\caption{(Color online) Pictorial illustration of the random variables used in the text. The light blue dots represent the observation points $t_1$, $t_2$ and $t_3$.
The red squares are the epochs $T_0 = 0, T_1, \ldots, T_5$. The conditional residual life-time is the time elapsed between 
$t_i$ and the next epoch $T_{n_i+1}$. It depends on previous values of $n_i$, this is the number of events between $0$ and $t_i$, 
with the event at $t= T_0 = 0$ not considered. Here, we have $n_1 = 1$,
$n_2 = 2$ and $n_3$ = 4. All the equations in this paper can be derived by analyzing this figure.}
\label{fig:XX}
\end{figure}

%Characterization of the fractional Poisson process
As a consequence of Kolmogorov's extension theorem, in order to fully characterize the stochastic process $\Nb (t)$, one has to 
derive its finite dimensional distributions. A further requirement on the process' paths
uniquely determines the process, namely that they are right-continuous step functions with left limits \cite{Billingsley1986}.
The finite-dimensional distributions are the multivariate probability distribution functions 
$\P (\Nb (t_1) = n_1, \Nb (t_2) = n_2, \ldots, \Nb(t_k) = n_k)$ with $t_1 < t_2 < \ldots < t_k$ and $n_1 \leq n_2 \leq \ldots \leq n_k$.
We have already given the formula for the one-point functions in Eq.~\eqref{eq:fractionalpoissondistribution}. 
The general finite dimensional distribution can be computed observing that the event
$\{\Nb (t_1) = n_1, \Nb (t_2) = n_2, \ldots, \Nb(t_k) = n_k\}$ is equivalent to
$\{0< T_{n_1} < t_1, T_{n_1+1} > t_1, t_1 < T_{n_2} < t_2, T_{n_2 + 1} > t_2, \ldots, t_{k-1} < T_{n_k} < t_k, T_{n_k+1}>t_k \}$.
Therefore, we find
\begin{multline}
\label{generalfidi}
\P(\Nb (t_1) = n_1, \Nb (t_2) = n_2, \ldots, \Nb(t_k) = n_k) = \\
\P(0< T_{n_1} < t_1, T_{n_1+1} > t_1, t_1 < T_{n_2} < t_2, T_{n_2 + 1} > t_2, \\
\ldots,  t_{k-1} < T_{n_k} < t_k, T_{n_k+1}>t_k) = \\
 \int_0^{t_1} du_1 f_\tau^{*n_1} (u_1) \int _{t_1 - u_1}^\infty du_2 f_\tau (u_2) \\
\int_{t_1 - u_1 - u_2}^{t_2 - u_1 - u_2} du_3 
f_\tau^{*(n_2-n_1 -1)}(u_3) \int_{t_2 - u_1 - u_2 - u_3}^\infty du_4 f_\tau (u_4)  \\
\ldots \int_{t_{k-1} - \sum_{i=1}^{2k-2} u_i}^{t_k - \sum_{i=1}^{2k-2} u_i} du_{2k-1} f_\tau^{*(n_k - n_{k-1} -1)} (u_{2k-1})  \\
\left[ 1 - F_\tau \left( t_k - \sum_{i=1}^{2k-1} u_i \right) \right].
\end{multline}
For instance, the two point function is given by
\begin{multline}
\label{twopointfidi}
\P(\Nb(t_1) = n_1, \Nb(t_2) =n_2) = \\
\P(0< T_{n_1} < t_1, T_{n_1+1} > t_1, t_1 < T_{n_2} < t_2, T_{n_2 + 1} > t_2) = \\
 \int_0^{t_1} du_1 f_\tau^{*n_1} (u_1) \int _{t_1 - u_1}^\infty du_2 f_\tau (u_2) \\
\int_{t_1 - u_1 - u_2}^{t_2 - u_1 - u_2} du_3 
f_\tau^{*(n_2-n_1 -1)}(u_3) \\
\left[ 1 - F_\tau \left( t_2 - u_1 - u_2 - u_3 \right) \right].
\end{multline}

%Increments and predictive probabilities
Let us focus on the two-point case for the sake of illustration.
As $\Nb(t)$ is a counting process, one has $\mathbb{P}(\Nb (t_1) = n_1, \Nb (t_2) = n_2) =
\P(\Nb (t_1) = n_1, \Nb (t_2) - \Nb (t_1) = n_2 - n_1)$ and, as a consequence of the
definition of conditional probability
\begin{multline}
\label{eq:twopoints}
\P(\Nb (t_1) = n_1, \Nb (t_2) - \Nb (t_1) = n_2 - n_1) = \\
\P(\Nb (t_2) - \Nb (t_1) = n_2 - n_1|\Nb (t_1) = n_1) \times \\
\times \P (\Nb (t_1) = n_1).
\end{multline}
For $\beta =1$, when the fractional Poisson process coincides with the standard Poisson process, the increments are
iid random variables and one has
\begin{multline}
\label{poissonprocessfidi}
\P(N_1 (t_2) - N_1 (t_1) = n_2 - n_1|N_1 (t_1) = n_1) = \\
\P(N_1 (t_2) - N_1 (t_1) = n_2 - n_1) = \\
\exp(-(t_2-t_1)) \frac{(t_2 - t_1)^{(n_2-n_1)}}{(n_2 - n_1)!}.
\end{multline}
On the contrary, for $0 < \beta < 1$, the increment $\Nb (t_2) - \Nb (t_1)$ and 
$\Nb (t_1)$ are not independent. Note that $\Nb(t_1)$ can be seen
as an increment as $\Nb(0) = 0$ by definition. However from Eq. \eqref{eq:twopoints}, the conditional probability
of having $n_2 -n_1$ epochs in the interval $(t_1,t_2)$ conditional on the observation of $n_1$ epochs in the interval
$(0,t_1)$ can be written
as a ratio of two finite dimensional distribution:
\begin{multline}
\label{eq:increment}
\P(\Nb (t_2) - \Nb (t_1) = n_2 - n_1|\Nb (t_1) = n_1) = \\
\frac{\P(\Nb(t_1)=n_1,\Nb(t_2)=n_2)}{\P (\Nb (t_1) = n_1)}.
\end{multline}
This probability can be evaluated by means of an alternative method, more appealing for a direct and 
practical understanding of the dependence structure. Let 
\bq
Y_{n_1} \stackrel{\mathrm{def}}{=} [T_{n_1 +1} - t_1 | \Nb(t_1)=n_1]
\eq
denote the residual lifetime at time $t_1$ (that is the time to the next epoch or renewal) conditional
on $\Nb (t_1) = n_1$ . With reference to Fig. \ref{fig:XX}, one can see that the conditional
probability $\P(\Nb (t_2) - \Nb (t_1) = n_2 - n_1|\Nb (t_1) = n_1)$ is given by
the following convolution integral for $n_2 - n_1 \geq 1$
\begin{multline}
\P(\Nb(t_{2})-\Nb(t_{1})=n_{2}-n_{1}|\Nb(t_{1})=n_{1})=\\
\int_{0}^{t_{2}-t_{1}} \mathbb{P}(\Nb(t_{2}-t_{1}-y)=n_{2}-n_{1}-1) 
f_{Y_{n_1}} (y) \, dy,
\label{eq:P2}
\end{multline}
where $f_{Y_{n_1}}(t)$ is the pdf of $Y_{n_1}$. In the case $n_2 - n_1 = 0$, one has
\bq
\P(\Nb(t_{2})-\Nb(t_{1})=0|\Nb(t_{1})=n_{1}) = 1-F_{Y_{n_1}}(t_2-t_1)
\eq
where $F_{Y_{n_1}} (y)$ is the cdf of $Y_{n_1}$. The distribution of the conditional residual lifetime $Y_{n_1}$ 
can be evaluated in several ways. For instance, one can notice that it can be decomposed as follows
\bq
\label{ydecomposition}
Y_{n_1} = \tilde\tau_{n_1+1}+U_{n_1} 
\eq
where $U_{n_1}$ is defined as
\bq
U_{n_1} \stackrel{\mathrm{def}}{=} [T_{n_1} |\Nb(t_1) = n_1],
\eq
and is the position of the last epoch before $t_1$ conditional on $\Nb(t_1)=n_1$, and 
\bq
\tilde\tau_{n_1+1} \stackrel{\mathrm{def}}{=}  [\tau_{n_1+1} -t_1 |T_{n_1 +1}>t_1]
\eq
is the difference between $\tau_{n_1+1}$ and $t_1$ conditional on $T_{n_1 + 1}>t_1$. 
The pdf of $U_{n_1}$ is given by the following chain of equalities
\begin{multline}
\label{U1}
\begin{aligned}
f_{U_{n_1}}(t) dt = & \P(T_{n_1} \in dt| \Nb(t_1) = n_1 )\\
= &  \P(T_{n_1} \in dt|T_{n_1}<t_1,T_{n_1}+\tau_{n_1 +1}> t_1)\\
= & \P(T_{n_1} \in dt|T_{n_1}<t_1,\tau_{n_1 +1}> t_1- T_{n_1})\\
\stackrel{\star}{=} & \frac{ \P(T_{n_1}\in dt)\int_{t_1-t}^\infty \P(\tau_{n_1+1}\in dw)}{ \P(T_{n_1}<t_1,\tau_{n_1 +1}>t_1 - T_{n_1})}\\
\stackrel{\ast}{=} & \frac{f_\tau^{*n_1}(t) [1-F_\tau(t_1 - t)] dt}{\int_0^{t_1}du  f_\tau^{*n_1}(u) [1-F_\tau(t_1 - u)]},
\end{aligned}
\end{multline}
where we used the independence between $T_{n_1}$ and $\tau_{n_1+1}$ $(\star)$ and  $f_{T_{n_1}}(x) = f_{\tau}^{* n_1}(x)$ $(\ast)$. 
The pdf of $\tilde{\tau}_{n_1 +1}$ is
\begin{multline}
\label{eq:tautilde}
\begin{aligned}
f_{\tilde\tau_{n_1+1}} (t|U_{n_1}) dt = & \P(\tau_{n_1+1} - t_1 \in dt  |T_{n_1+1}>t_1)\\
=& \frac{ \P(\tau_{n_1+1} \in dt + t_1)}{\P(\tau_{n_1+1}>t_1 - U_{n_1})} \\
=&\frac{f_\tau(t+t_1) dt}{1-F_\tau(t_1- U_{n_1})}.
\end{aligned}
\end{multline}
From Eq. \eqref{ydecomposition}, one can write that
\bq
f_{Y_{n_1}} (t) = \int_{0}^{t_1} f_{\tilde\tau_{n_1+1}}(t-u|u)f_{U_{n_1}}(u) du
\eq
and this equation leads to
\bq
f_{Y_{n_1}} (t) = 
\frac{\int_{0}^{t_1} \, du f_{\tau}^{*n_1} (u) f_\tau (t+t_1-u)}{\int_{0}^{t_1} \, du f_{\tau}^{*n_1} (u) [1 - F_\tau (t_1 - u)]}
\label{eq:fYn1final}
\eq
that, together with Eq.~\eqref{eq:fractionalpoissondistribution}, gives us the probability of the conditional
increments in Eq.~\eqref{eq:P2}. Notice that, for $n_1 = 0$, one has $f_\tau^{*0} (u) = \delta(u)$ and Eq.~\eqref{eq:fYn1final} 
reduces to the familiar equation for the residual life-time pdf in the absence of previous renewals
\begin{equation}
f_{Y_0}(t) = \frac{f_\tau (t+t_1)}{1 - F_\tau (t_1)}.
\label{eq:fYn1zero}
\end{equation}

This method can be applied to the general multidimensional case. As in Eq.~\eqref{eq:twopoints} we can write
\begin{multline}
\P(\Nb(t_1)=n_1,\ldots,\Nb(t_k)=n_k,\Nb(t_{k+1})=n_{k+1}) =\\
\P(\Nb(t_{k+1})-\Nb(t_{k})=n_{k+1}-n_k| \\
\Nb(t_1)=n_1,\ldots,\Nb(t_k)=n_k)\times\\
\times \P(\Nb(t_1)=n_1,\ldots,\Nb(t_k)=n_k)
\label{eq:multifidi}
\end{multline}
and the predictive probabilities can be evaluated as
\begin{multline}
\P(\Nb(t_{k+1})-\Nb(t_{k})=n_{k+1}-n_k|\ldots\\
|\Nb(t_1)=n_1,\ldots,\Nb(t_k)=n_k) =\\
\int_0^{t_{k+1}-t_k} \P(\Nb(t_{k+1}-t_k-y)=n_{k+1}-n_k-1) \times \\
\times f_{Y_{n_1,\ldots,n_k}}(y)dy,
\label{eq:multifidicond}
\end{multline}
where we defined
\bq
Y_{n_1,\ldots,n_k} \stackrel{\mathrm{def}}{=} [T_{n_k +1}-t_k | \Nb(t_1)=n_1, \ldots, \Nb(t_k)=n_k].
\eq
Again, we can use a decomposition of $Y_{n_1,\ldots,n_k}$
\bq
Y_{n_1,\ldots,n_k} = \tilde\tau_{n_k+1} + U_{n_k},
\eq
where
\bq
U_{n_k} \stackrel{\mathrm{def}}{=} [T_{n_k} |  \Nb(t_1)=n_1, \ldots, \Nb(t_k)=n_k],
\eq
and
\bq
\tilde\tau_{n_k+1} \stackrel{\mathrm{def}}{=} [\tau_{n_k+1} - t_k| T_{n_k+1}>t_k].
\eq
The difference with the two-point case is that $U_{n_1} = [T_{n_1}|\Nb(t_1)=n_1] = 
\left[\sum_{i=1}^{n_1} \tau_i | \Nb(t_1)=n_1\right]$ must be replaced by
\bq
U_{n_k} = t_{k-1}+Y_{n_1,\ldots,n_{k-1}} + \left[ \sum_{i=n_{k-1}+1}^{n_k} \tau_i  |  \Nb(t_k)=n_k \right].
\eq
The time between $t_{k-1}$ and the next renewal epoch is $Y_{n_1,\ldots,n_{k-1}}$ and it is independent from $\sum_{i=n_{k-1}+1}^{n_k} \tau_i$. Therefore,
the convolution 
\bq
q(n_1,\ldots,n_k;t) = f_{Y_{n_1,\ldots,n_{k-1}}}*f_{\tau}^{*(n_k-n_{k-1}-1)}(t)
\eq
replaces $f_\tau^{*n_1}(t)$ in Eq. \eqref{U1}.
This leads to 
\begin{multline}
f_{U_{n_k}}(z) =\\
 \frac{q(n_1,\ldots,n_k;t+t_{k-1})[1-F_\tau (t_k-t)]}{\int_{t_{k-1}}^{t_k} q(n_1,\ldots,n_k;u+t_{k-1})[1-F_\tau (t_k-u)] du  }.\\
\end{multline}
On the other hand, $f_{\tilde\tau_{n_k+1}}(t)$ has the same functional form as $f_{\tilde\tau_{n_1+1}}(t)$ given in Eq.~\eqref{eq:tautilde} with
$U_{n_k}$ replacing $U_{n_1}$. Therefore, $Y_{n_1,\ldots,n_k}$ has the following pdf
\begin{multline}
f_{Y_{n_1,\ldots,n_k}}(t) = \\
\frac{\int_{t_{k-1}}^{t_k} du \, q(n_1,\ldots,n_k;u+t_{k-1}) f_\tau (t+t_{k} - u)}
{\int_{t_{k-1}}^{t_k} du \, q(n_1,\ldots,n_k;u+t_{k-1})[1-F_\tau (t+t_k - u)]}.
\label{eq:Yk}
\end{multline}
In practice, the random variable $Y_{n_1,\ldots,n_{k-1}}$ carries the memory of the observations made at times 
$t_1,\ldots,t_{k-1}$; the knowledge of $f_{Y_{n_1,\ldots,n_{k-1}}}$ allows the computation of $f_{Y_{n_1,\ldots,n_{k}}}$, and, 
via Eqs.~\eqref{eq:multifidi} and \eqref{eq:multifidicond}, the $k+1$-dimensional distribution can be derived as well.
 
\begin{figure*}[p]
\includegraphics[width=0.7\columnwidth]{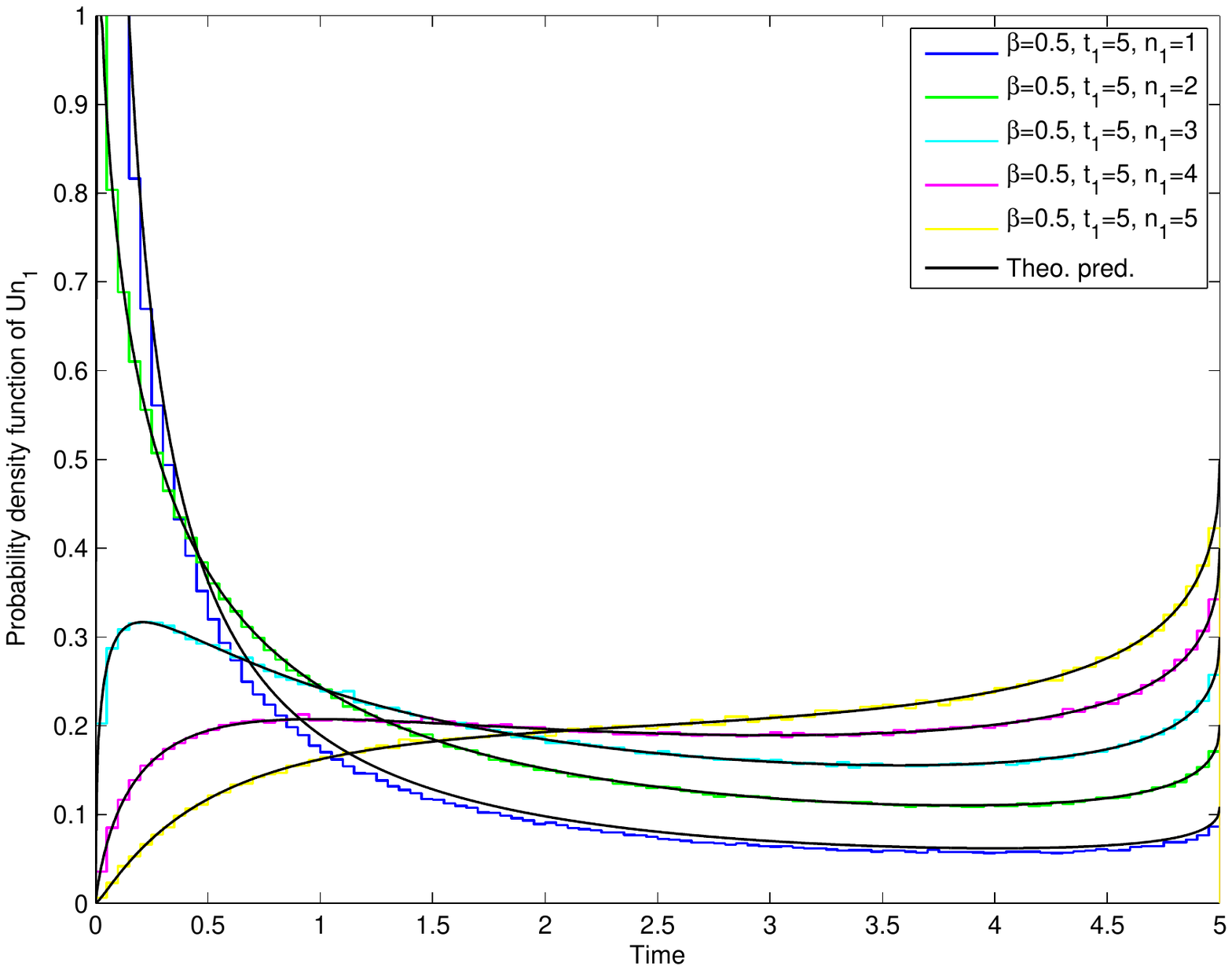}
\includegraphics[width=0.7\columnwidth]{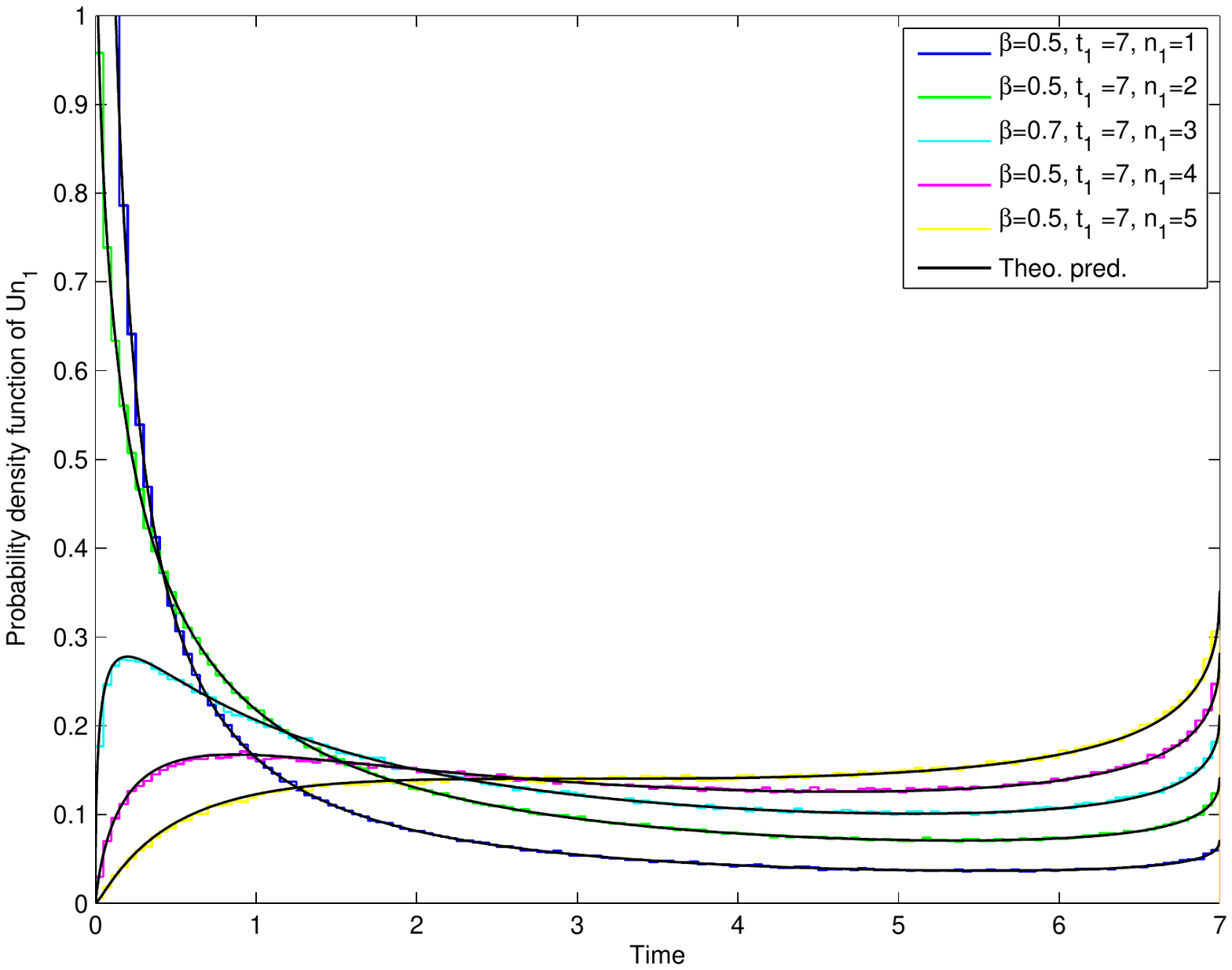}
\includegraphics[width=0.7\columnwidth]{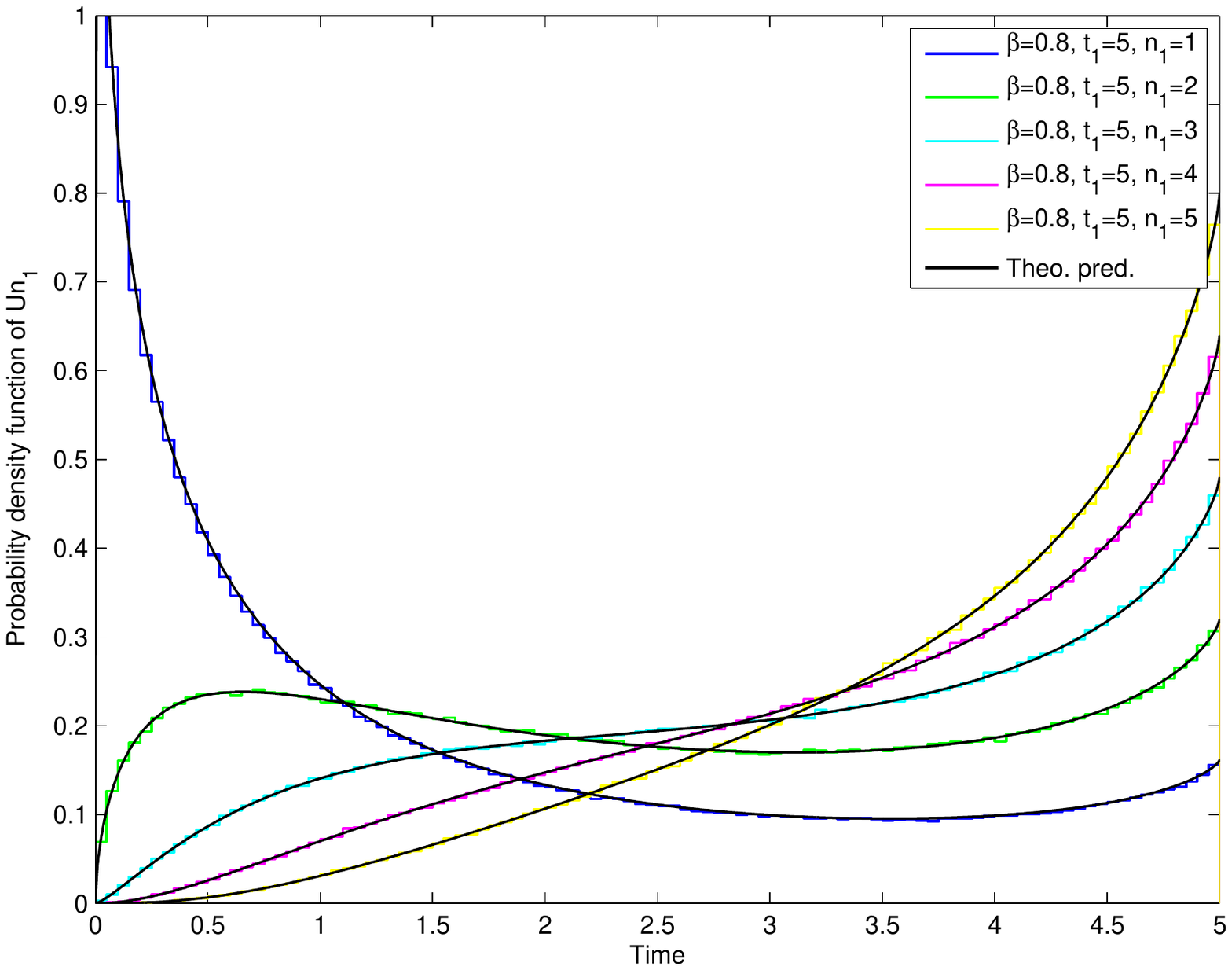}
\includegraphics[width=0.7\columnwidth]{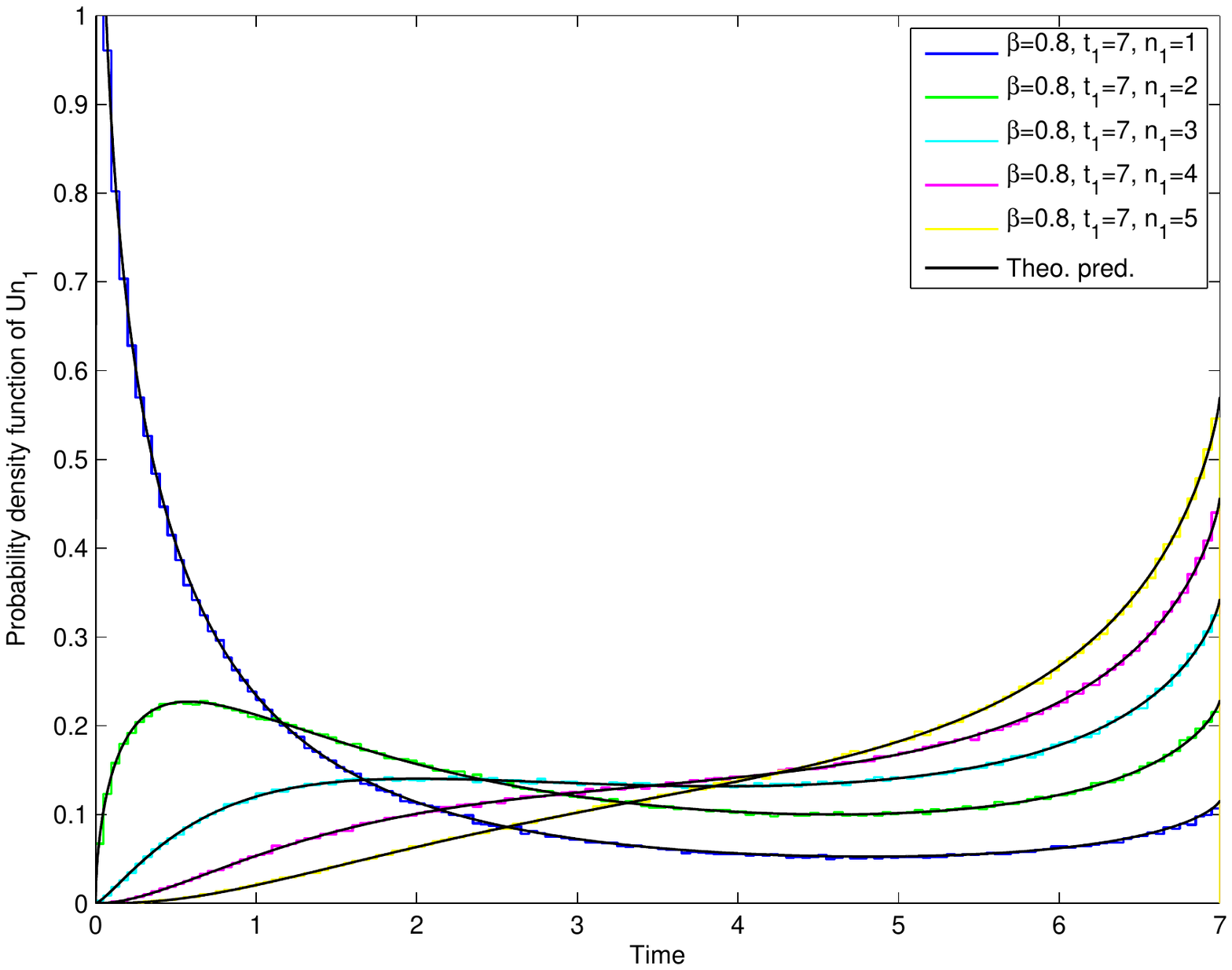}
\includegraphics[width=0.7\columnwidth]{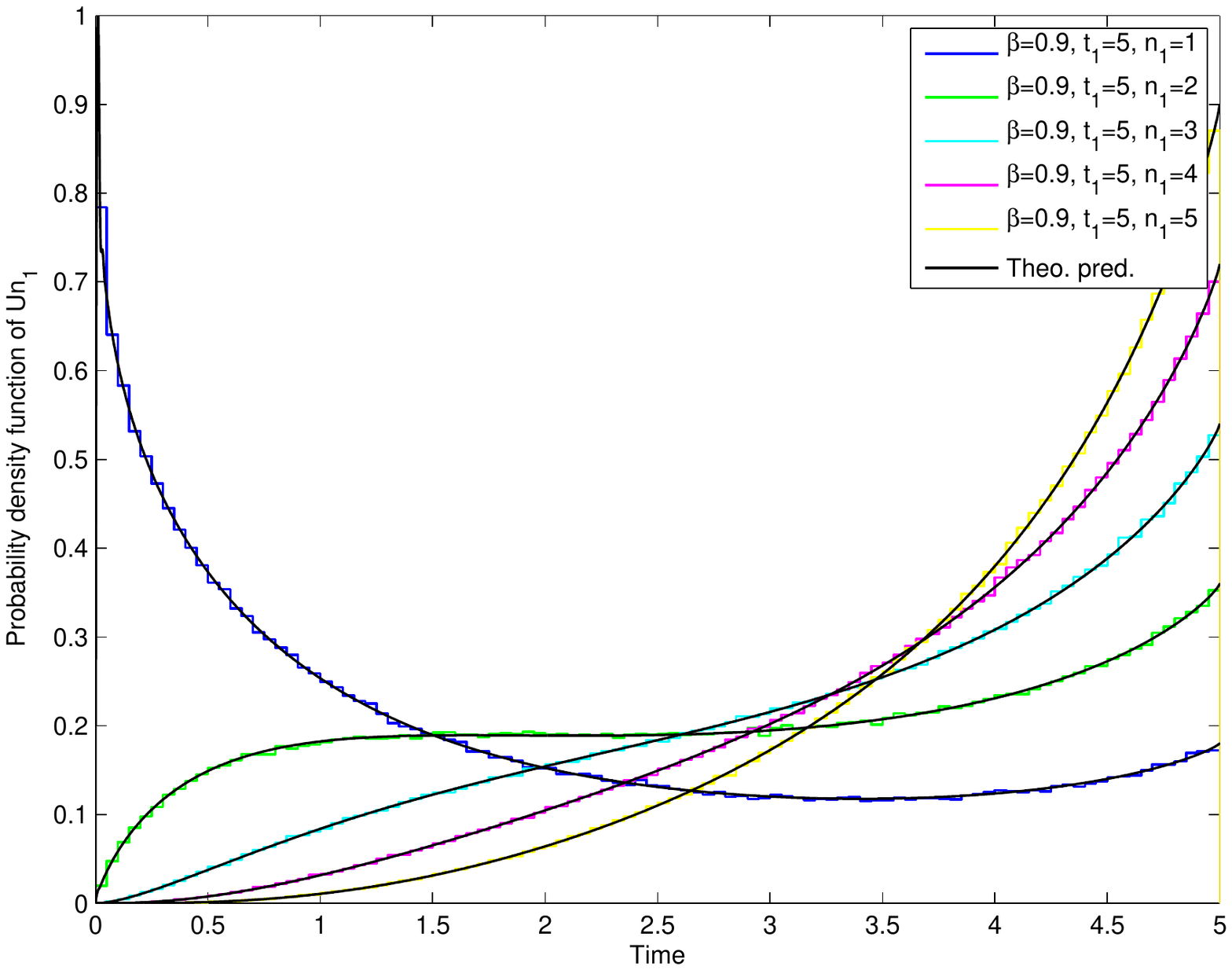}
\includegraphics[width=0.7\columnwidth]{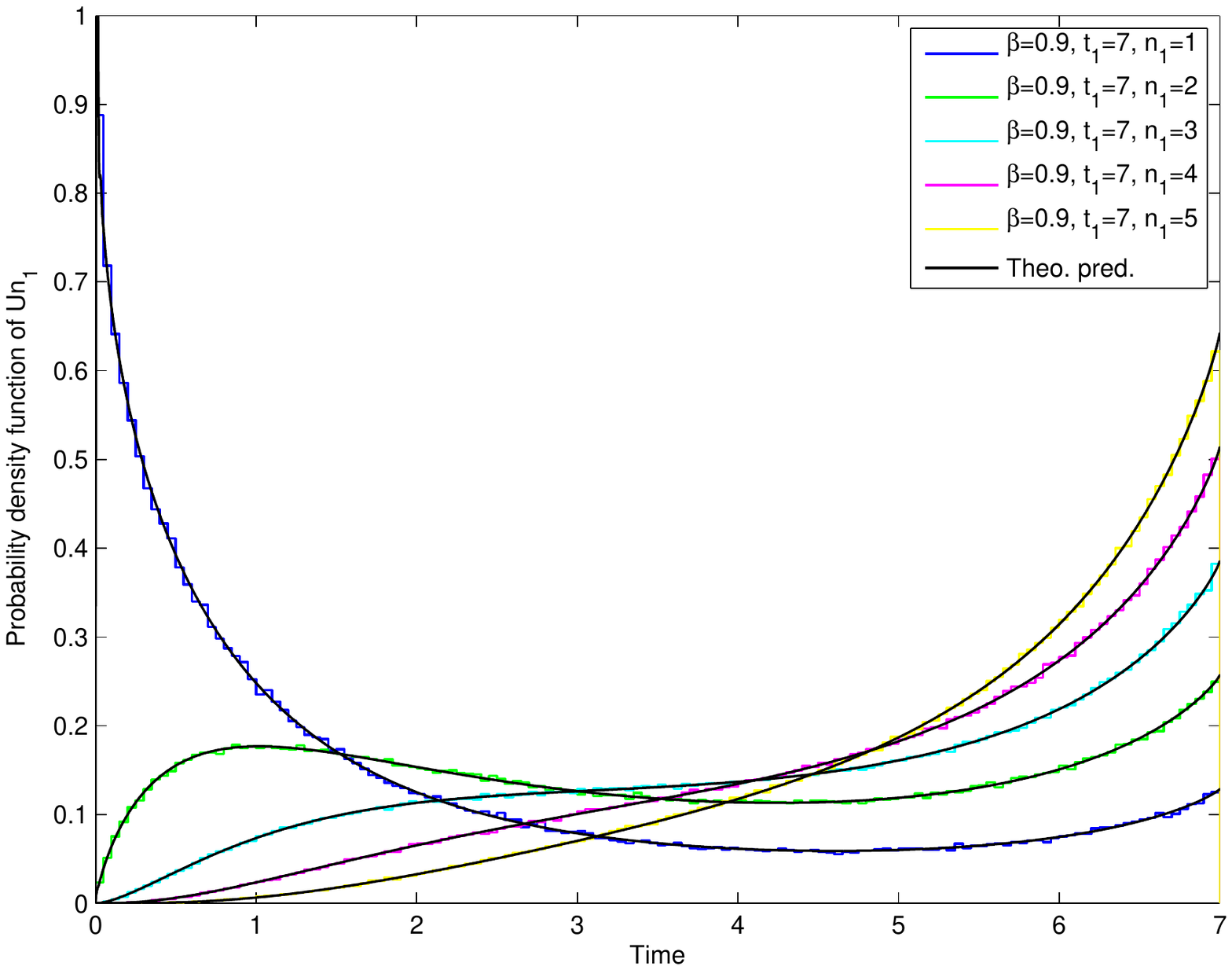}
\caption{(Color online) Pdf of the random variable $U_{n_1}$ as given in Eq.~\eqref{U1} (solid black lines) compared 
to Monte Carlo simulations (colored step lines) for three values of $\beta$ and two different values of $t_1$. 
$10^7$ different paths were simulated for each value of $\beta$ and the bin width is 0.05. Time is in arbitrary units. }
\label{fig:prob1U}
\end{figure*}

\begin{figure*}[p]
\includegraphics[width=0.7\columnwidth]{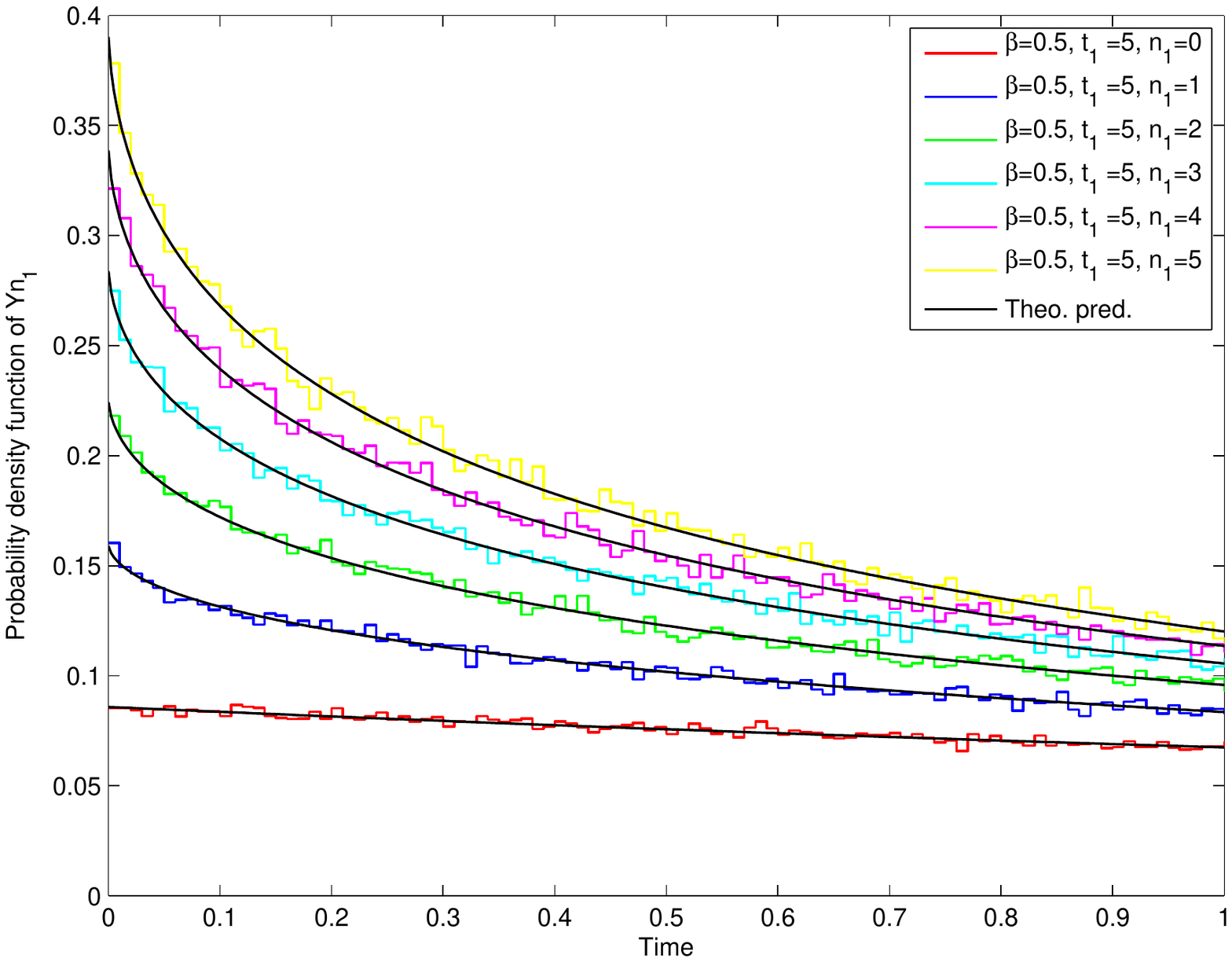}
\includegraphics[width=0.7\columnwidth]{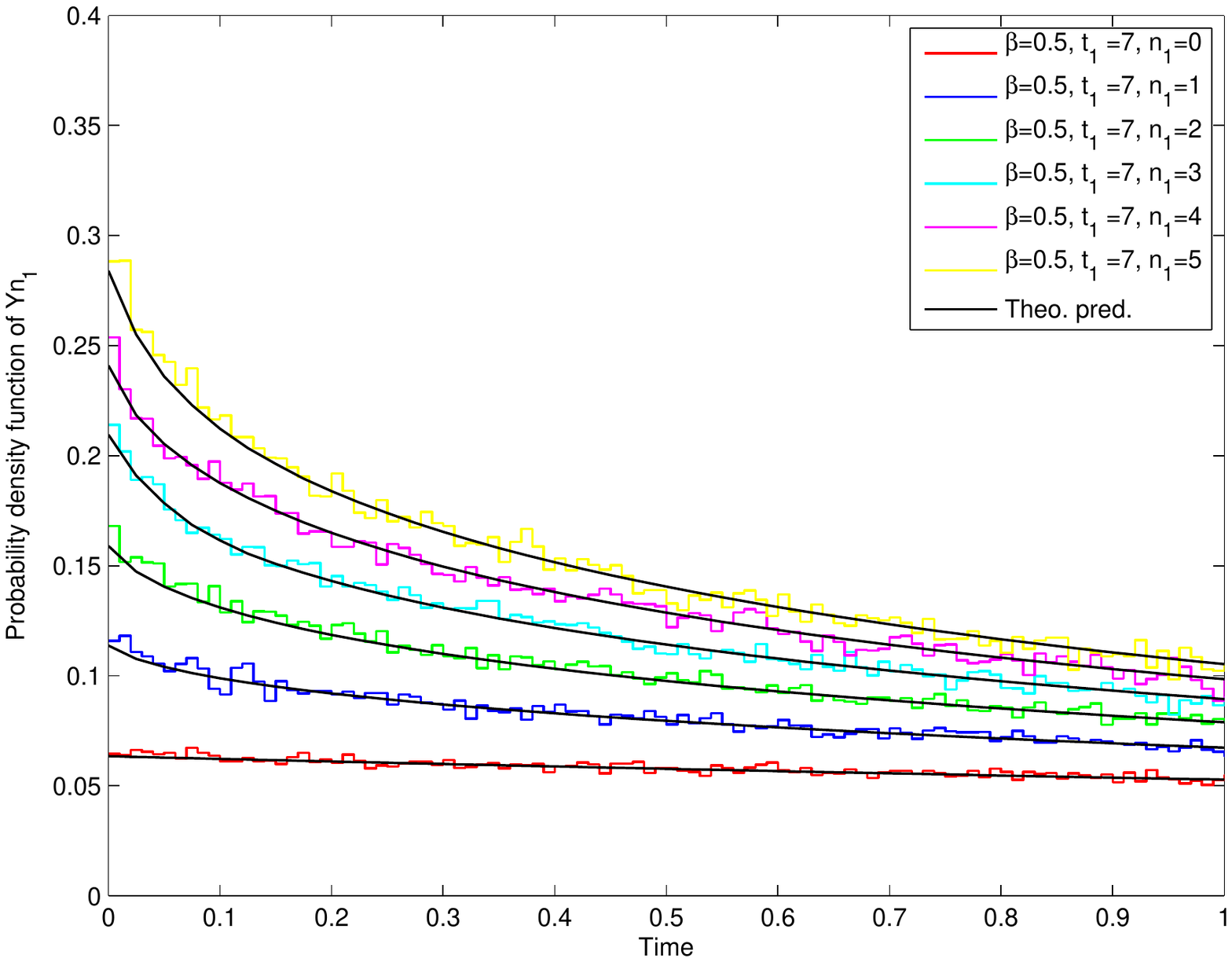}
\includegraphics[width=0.7\columnwidth]{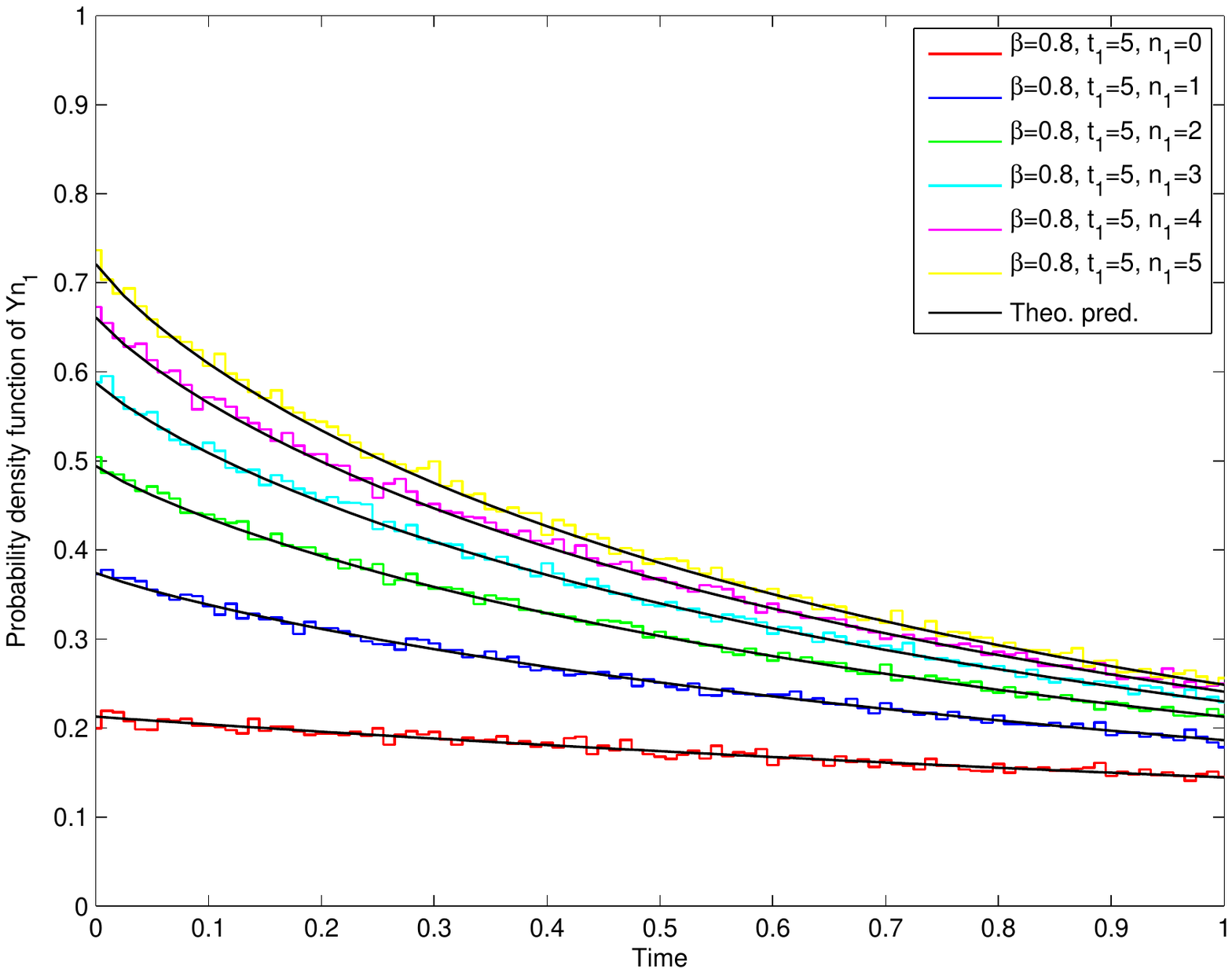}
\includegraphics[width=0.7\columnwidth]{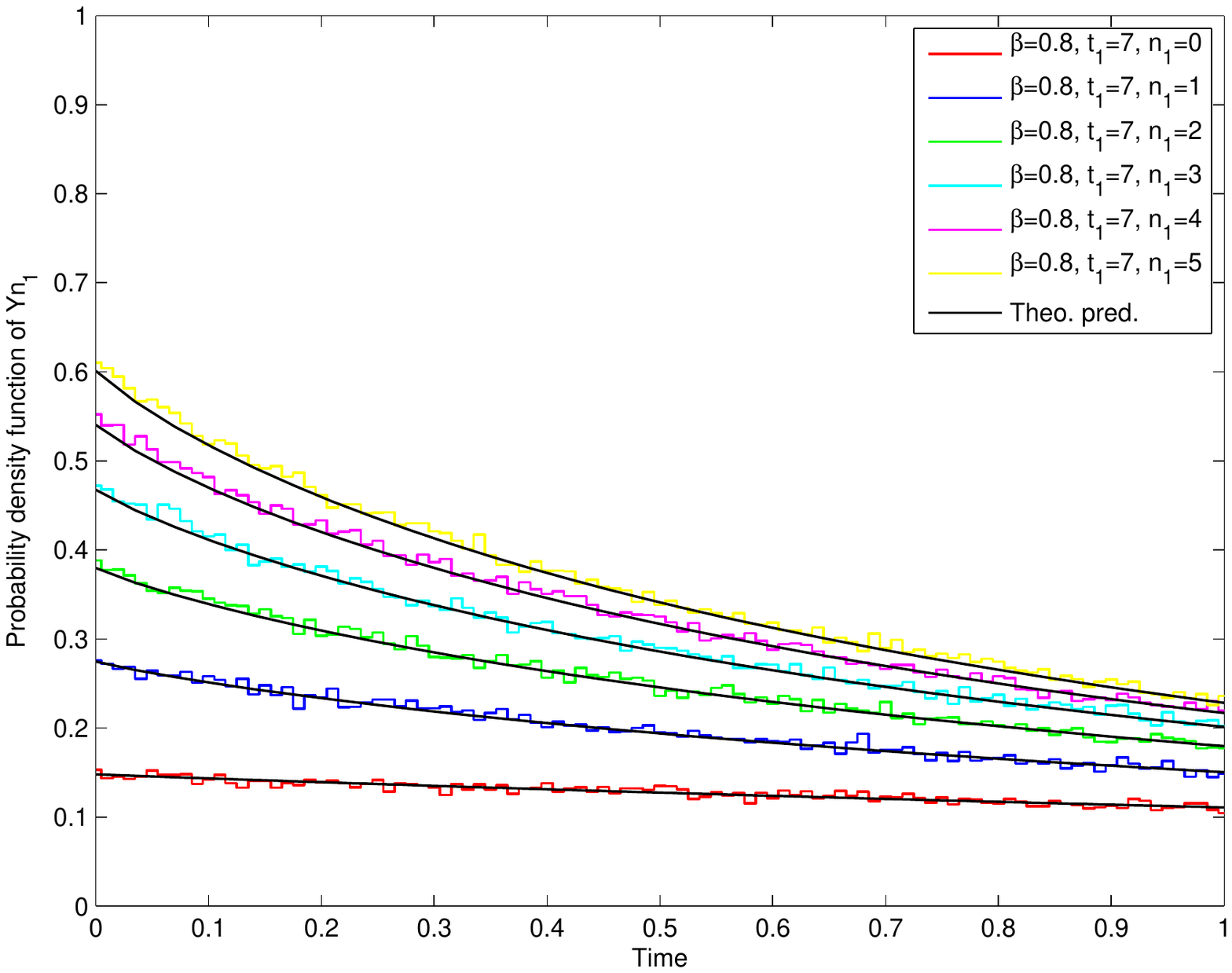}
\includegraphics[width=0.7\columnwidth]{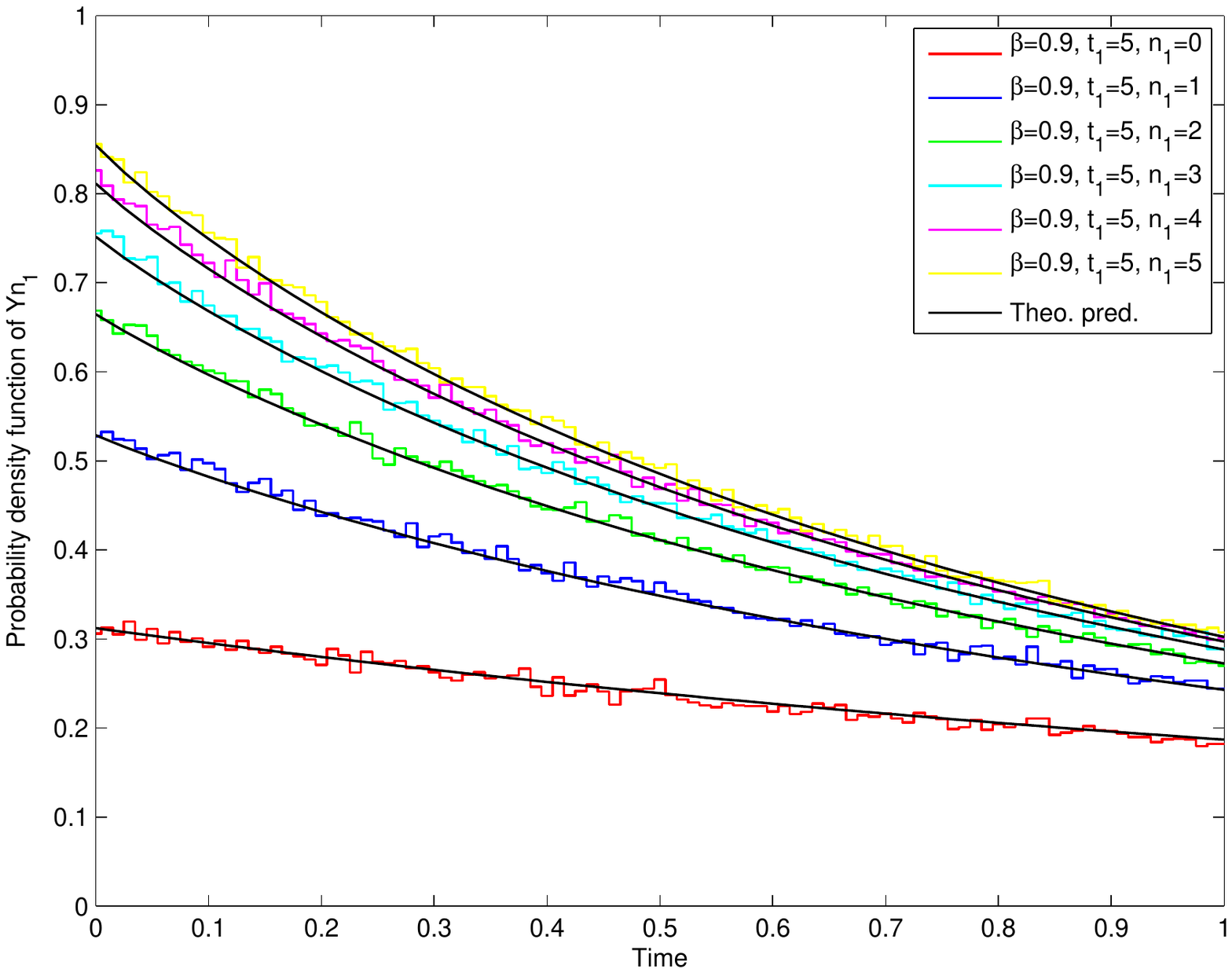}
\includegraphics[width=0.7\columnwidth]{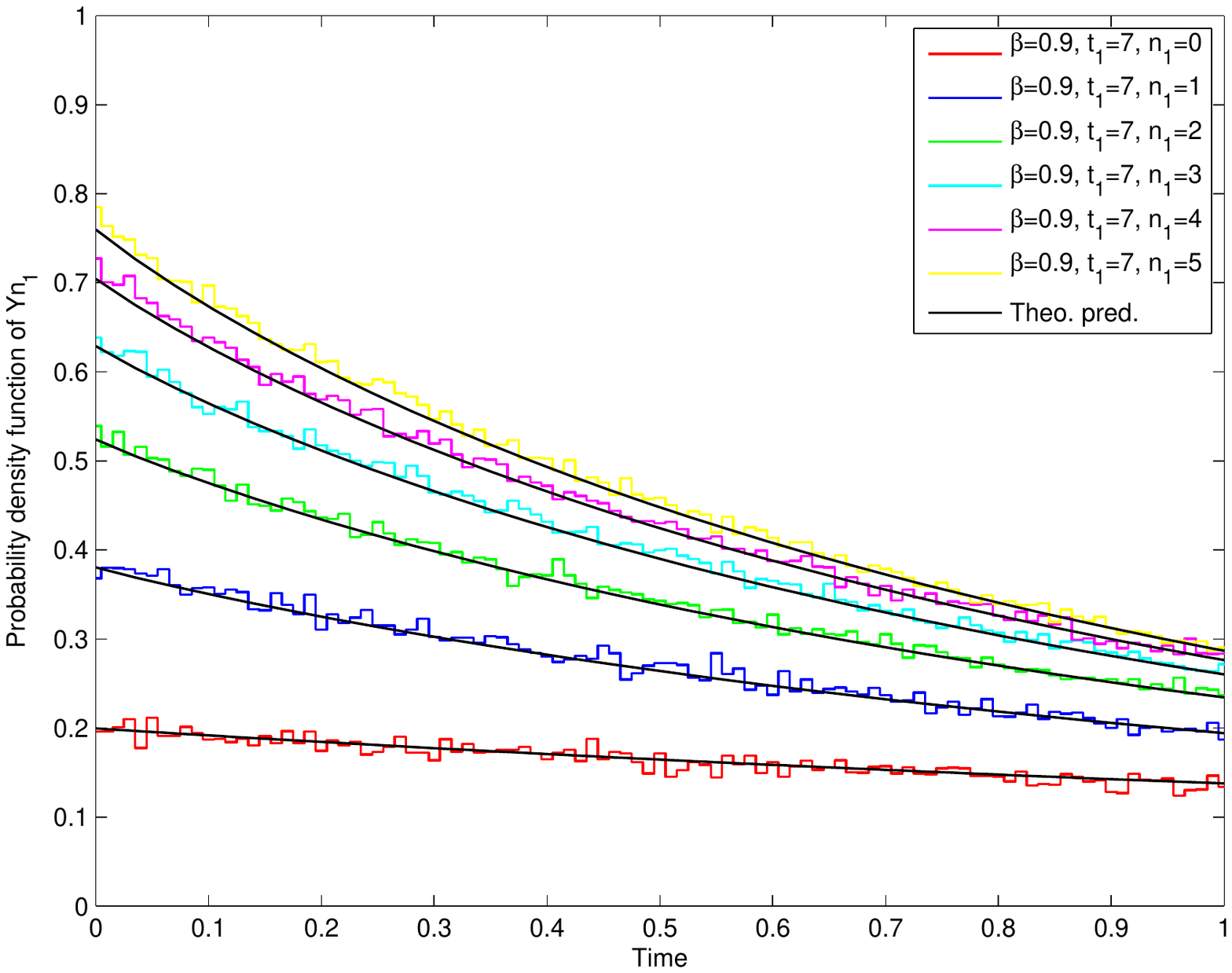}
\caption{(Color online) Pdf of the random variable $Y_{n_1}$ as given in Eqs.~\eqref{eq:fYn1final} and \eqref{eq:fYn1zero} (solid black lines)
compared to Monte Carlo simulations (colored step lines) for three values of $\beta$ and two different values of $t_1$. $10^7$ 
different paths were simulated for each value of $\beta$ and the bin width is 0.01. Time is in arbitrary units.}
\label{fig:prob1Y}
\end{figure*}

Figs.~\ref{fig:prob1U} and \ref{fig:prob1Y} compare the theoretical results of Eqs.~\eqref{U1}, \eqref{eq:fYn1final} 
and \eqref{eq:fYn1zero} with those of a Monte Carlo simulation based on the algorithm presented in equation (20) of Ref. \cite{Fulger2008}. 

%Acks
The Japanese Society for the Promotion of Science (grant N. PE09043) supported MP during his stay at the International Christian University
in Tokyo, Japan.

\bibliographystyle{apsrev4-1}
\bibliography{fractionalPoisson} 
\end{document}